\documentclass[a4paper,pdftex,Skim]{amsart}

\usepackage[german]{babel}
\usepackage[fixlanguage]{babelbib}
\selectbiblanguage{german}
\usepackage{geometry}
\usepackage[utf8]{inputenc}
\usepackage{textcomp}
\usepackage[T1]{fontenc}
\usepackage{afterpage}
\usepackage[numbers,square]{natbib}
\usepackage[scaled=0.9]{helvet}
\usepackage{mathptmx}
\usepackage{courier}
\usepackage{microtype}
\usepackage[font=footnotesize, normal, bf, up, justification=raggedright]{caption}[2008/08/24]
\usepackage{graphicx}
\usepackage{hyperref}
\hypersetup{
	colorlinks=true,
	pdfauthor={M. Kiessling, T. Kreisel, S. Kurz, J. Rambau, K. Schade, C. Schwarz},
	pdftitle={Das Optimierungslabor -- ein Erfahrungsbericht},
}

\begin{document}

\title{Das Optimierungslabor -- ein Erfahrungsbericht}

\author{Miriam Kie\ss ling, Sascha Kurz, Tobias Kreisel, J\"org Rambau, Konra Schade, und Cornelius Schwarz}

\begin{abstract}
For several years, students visit us on different occasions at the university. But how to bridge from the
school curriculum to the contents of the university mathematics? And how to find a focal point at which an 
active contribute, despite the lack of knowledge, in view of limited time is possible? Our approach: Translate, 
under guidance, everyday life optimization problems into the language of
mathematics, i.e. using variables, target functions, equations and inequalities. These so-called
integer linear programming models are then solved by standard software. In this report we wnat to tell about
the lessons we have learned.

\bigskip

Seit mehreren Jahren besuchen uns Sch\"ulerinnen und Sch\"uler zu
unterschiedlichen Anl\"assen an der Universit\"at.  Doch wie l\"asst sich die
Br\"ucke vom Schulstoff zu den Inhalten der Universit\"atsmathematik schlagen?
Und: findet man einen Themenschwerpunkt, bei dem ein aktives Mitmachen trotz
fehlender Vorkenntnisse in Anbetracht begrenzter Zeit m\"oglich wird? Unser
Ansatz: Unter Anleitung sollen Alltagsoptimierungsprobleme in die Sprache der
Mathematik \"ubersetzt werden, also durch Variablen, Bewertungsfunktionen,
Gleichungen und Ungleichungen ausgedr\"uckt werden.  Diese sogenannten
ILP-Modelle werden dann durch Standardsoftware der diskreten Optimierung
gel\"ost und deren L\"osung wieder in die Alltagssituation zur\"uck\"ubersetzt.
In diesem Bericht wollen wir unsere Erfahrungen mit konkreten Details der
Umsetzung schildern.
\end{abstract}

\maketitle

\section{Einleitung}
\label{sec:introduction}

Seit mehreren Jahren besuchen uns Sch\"ulerinnen und Sch\"uler der
Jahrgangsstufen 10--13 an der Universit\"at Bayreuth zu Anl\"assen 
wie dem Tag der
Mathematik\footnote{\,\url{http://www.tdm.uni-bayreuth.de}},
dem Girls' Day, der MINT-Herbstuniversit\"at oder einfach auf Initiative ihrer 
Klassenleitungen.  Sie m\"ochten einen Einblick in die Welt der
Mathematik \"uber die Schulmathematik hinaus bekommen. Doch wie l\"asst
sich die Br\"ucke vom Schulstoff zu den Inhalten der
Universit\"atsmathematik schlagen? Und: findet man einen
Themenschwerpunkt, bei dem ein aktives Mitmachen trotz fehlender
Vorkenntnisse und in Anbetracht begrenzter Zeit m\"oglich wird?

In der diskreten Optimierung lassen sich Problem-Modellierung und
Problem-L\"osung sehr gut trennen. Selbst forschungsnahe Modelle der
ganzzahligen linearen Optimierung (ILP-Modelle) basieren auf sehr
elementaren \"Uberlegungen, wie die Entscheidungsm\"oglichkeiten,
Ziele und Restriktionen eines Alltagsproblems in Variablen,
Bewertungsfunktionen, Gleichungen und Ungleichungen ausgedr\"uckt
werden k\"onnen. Wie dann optimale L\"osungen gefunden werden,
erfordert zwar tiefergehende Mathematik, es gibt aber Software daf\"ur,
in der das Wissen aus Teilen des Mathematik-Studiums und der
mathematischen Forschung kondensiert vorliegt.

Unser Vermittlungsziel: Sch\"ulerinnen und Sch\"uler wissen nach dem
Besuch, dass man verschiedenste Probleme angreifen kann, indem man sie
in die Sprache der Mathematik \"ubersetzt, denn in Software gegossenes
mathematisches Know-How kann dann diese Probleme l\"osen, ohne etwas
\"uber die Probleme selbst zu wissen. Unsere Idee f\"ur eine Ma\ss nahme:
Ein Optimierungslabor. Die Sch\"ulerinnen und Sch\"uler isolieren in
Teamarbeit die wesentlichen logischen Merkmale von Sudokul\"osen,
Rucksackpacken, Routenplanung u.\,v.\,a.\,m. Dann \"ubersetzen sie die
Problemstellungen in die Sprache der Mathematik (hier: ILP-Modelle)
und lassen sie (unterst\"utzt durch unser Team) von Computerprogrammen
l\"osen (ILP-L\"oser), die nichts anderes als diese Sprache verstehen.
Schlie\ss lich \"ubersetzen sie die mathematischen L\"osungen wieder in die
Sprache der Problemstellung. Erfahrungen mit der Modellierung auf
Basis linearer Gleichungssysteme k\"onnen dabei aus dem Schulunterricht
eingebracht werden.

In diesem Bericht wollen wir unsere Erfahrungen mit konkreten Details
der Umsetzung schildern.

\section{Historie des Projekts}
\label{sec:hist-des-proj}

Die Idee, Modellierungswerkzeuge und Standardsoftware zu benutzen, um
ohne viel Vorlauf Studierenden in fr\"uhen Semestern den Zugang zur
Leistungsf\"ahigkeit einer mathematischen Methode zu ebnen, ist nicht
von uns.  In Arbeitsgruppen der diskreten Optimierung werden
mittlerweile an vielen Universit\"aten Modellierungsaufgaben regelm\"a\ss ig
in Vorlesungen verwendet.

Eine der Pionierveranstaltungen in diesem Bereich f\"ur Schulklassen,
die Modellierungswoche der TU~Kaiserslautern und der TU~Darmstadt in
Lambrecht, l\"auft seit 1993.\footnote{\,Siehe den \"Uberblick unter
\url{http://www.kfunigraz.ac.at/imawww/modellwoche/konzept4.html}.}
Die daten- und softwareorientierte Darstellung von Themen aus der
Kombinatorischen Optimierung anhand von Praxisprojekten f\"ur
Studierende ist eine langj\"ahrige Idee von Martin Gr\"otschel
gewesen, formuliert zuerst als ein Buchprojekt "`Combinatorial
Optimization at Work"', Projekt G1 aus der Gr\"undungsphase des Matheon
Berlin 2004, sp\"ater als ein Blockkursprojekt mit Online-Dokumentation.

Unser Start war ein Blockkurs f\"ur Studierende ab dem 2.~Semester Anfang 2005
an der Universit\"at Bayreuth.\footnote{\,\url{http://www.rambau.wm.uni-
bayreuth.de/Teaching/Uni_Bayreuth/WS_2004/Diskrete_Optimierung_Anwendungen/}}
Teile des \"Ubungsmaterials fanden auch Verwendung im ersten Blockkurs
"`Combinatorial Optimization at Work"' in Berlin, 2005.\footnote{\,\url{http
://co-at-work.zib.de/berlin/}}

Schnell wurde deutlich, dass die Vorgehensweise nach weiteren
Vereinfachungsschritten auch f\"ur Sch\"ulerinnen und Sch\"uler sowie f\"ur die
allgemeine \"Offentlichkeit interessant sein k\"onnte. Die zun\"achst angebotenen
Labore an den Bayreuther Tagen der Mathematik 2007 und 2008 waren noch recht
nah an die universit\"are Arbeitsweise angelehnt. Es war zu merken, dass ein
Ausschnitt aus einem Blockkurs allein noch kein ausreichendes Konzept f\"ur eine
zweist\"undige Veranstaltung sein kann. Im Jahr der Mathematik 2008 wurde dann
die Strategie im Rahmen einer Kooperation mit der Stadtbibliothek Bayreuth
\"uberarbeitet. Kernpunkte dieser \"Uberarbeitung waren die Bereitstellung von
m\"oglichst greifbaren Materialien zur besseren Aktivierung von kooperativem
Nachdenken, Dokumentieren, Verwerfen, Revidieren, Sichern.  Im Nachgang dazu
hat sich ein Portfolio \footnote{\,\url{http://www.wm.uni-
bayreuth.de/index.php?id=optlabor}} von beispielhaften Optimierungsprojekten
stabilisiert, in denen die \"Ubersetzung des (idealisierten) Anwendungsproblems
in die Sprache der Mathematik im Mittelpunkt steht.

Mittlerweile wird das Optimierungslabor, wann immer die Zeit es
zul\"asst, eingeleitet von einem separat f\"ur die
URANIA\footnote{\,\url{http://www.urania.de/die-urania/}}-Vortragsreihe
des \textsc{Matheon}\footnote{\,\url{http://www.matheon.de/}} Berlin
entstandenen Vortrag \"uber die Optimierung der {\glqq}Gelben Engel{\grqq} des
ADAC.\footnote{\,\url{http://www.wm.uni-bayreuth.de/fileadmin/Sonstiges/Folien/Schueler_ADAC_2004-11-16.pdf}}
Hier werden die Grundprinzipien mathematischen Modellierens am
Beispiel des Handlungsreisendenproblems detailliert vorgef\"uhrt und in
die Rahmenhandlung des ADAC-Praxisprojekts eingebettet.  Ziel ist die
Einstimmung und der erste Kontakt mit dem Thema "`Mathematik als
Sprache"'. Das Publikum wird dabei zum aktiven Mitdenken angeregt;
viele legen beim Vortrag ihre anf\"angliche Zur\"uckhaltung nach und
nach ab.

\section{Der Einf\"uhrungsvortrag: Gelbe Engel, ein Handlungsreisender und die Sprache der Mathematik}
\label{sec:Vortrag}
Zu Beginn der Veranstaltung wird die Gruppe auf das Thema Modellierung
eingestimmt. Alle Projekte, die die Gruppe sp\"ater aktiv bearbeiten
wird, sind akademische Modellierungsbeispiele.  Ziel ist ein Begreifen
von Prinzipien.  Echte, praxistaugliche Modelle sind komplizierter,
beziehen sich aber in Teilen immer wieder auf diese Prinzipien.  Der
Prozess, in dem man aus einem hoch-komplexen Anwendungsproblem der
realen Welt wichtige Modellierungsprinzipien extrahiert, wird im
Vortrag nachverfolgt.  Motto des Vortrags ist ein Zitat von Galilei
aus seiner {\glqq}Goldwaage{\grqq}.\footnote{\,\emph{Saggiatore}, zu finden unter
\url{http://www.liberliber.it/biblioteca/g/galilei/}} Sinngem\"a\ss 
behauptet Galilei, dass das
Universum nur verstanden werden kann, wenn man die Sprache der
Mathematik beherrscht.  Wir behaupten etwas weniger: Manche Probleme
der Praxis kann man nur dann fundiert l\"osen, wenn man sie vorher in
die Sprache der Mathematik \"ubersetzt - der Kern der mathematischen
Modellierung.

Aufh\"anger ist das Einsatzplanungsproblem der {\glqq}Gelben Engel{\grqq} des ADAC:
Wie wird entschieden, welches Hilfefahrzeug in welcher Reihenfolge
welchen Havaristen hilft, so dass einerseits jedem geholfen wird und
andererseits Kosten und Wartezeiten m\"oglichst gering sind?  Ein
Beispiel f\"ur ein \emph{Optimierungsproblem}.  In einer ersten
Argumentationslinie werden die \emph{Online-} und die
\emph{Echtzeitproblematik} herausgestellt: Man wei\ss  nicht, wer in der
Zukunft noch havarieren wird, und die Entscheidung \"uber die
Einsatzplanung muss noch w\"ahrend des Telefonats mit dem Havaristen
getroffen werden, damit eine gesch\"atzte Ankunftszeit mitgeteilt werden
kann. Als grunds\"atzliches Vorgehen wird die \emph{Reoptimierung}
erl\"autert: plane zu jedem Zeitpunkt optimal f\"ur den Fall, dass keine
neue Havaristen mehr anrufen.  Ruft doch einer an, dann wiederhole die
Planung.  Dies reduziert das Problem auf ein
\emph{Offline-Optimierungsproblem} auf Daten eines
\emph{Schnappschusses} des Systems.  An dieser Stelle braucht man ein
Computerprogramm, das f\"ur jeden Schnappschuss des Systems eine
optimale Einsatzplanung echtzeittauglich (d.\,h.{ }in etwa 10\,s)
berechnet.

An dieser Stelle wird der Anwendungskontext verlassen, weil man in
einem Modell f\"ur das ADAC-Problem den Wald vor lauter B\"aumen nicht
mehr s\"ahe (Genaueres findet sich in~\cite{Rambau+Schwarz:ADAC-Aufzuege-Modelle:2008}
und~\cite{Rambau:Gelbe-Engel:2010}). Stattdessen wird
das Problem durch gezieltes Weglassen einiger Aspekte auf ein
\emph{Handlungsreisendenproblem (TSP)} reduziert: Ein
Handlungsreisender sucht eine k\"urzeste, geschlossene Tour durch eine
Menge von St\"adten.  Und die Modellierung eines TSPs enth\"alt einige
sch\"one grunds\"atzliche Aspekte, deren genaue, m\"oglichst interaktiv
gefundene Formalisierung den Kern des Vortrags ausmacht.

Die Struktur des Modellierungsprozesses orientiert sich an folgenden
Schritten:
\begin{enumerate}
\item Wie stellt man die \emph{wesentlichen Aspekte des Problems} dar?
  Dies f\"uhrt auf \emph{gewichtete, ungerichtete Graphen} als zentrale
  mathematische Struktur.
\item Wie \emph{kodiert} man eine spezielle TSP-Tour? Dies f\"uhrt zu
  \emph{Adjazenzmatrizen}, die in einer Tabellenkalkulation gut
  visualisierbar sind.
\item Wie repr\"asentiert man eine \emph{unbekannte} TSP-Tour? Dies
  f\"uhrt auf Adjazenzmatrizen, die \emph{Variablen} als Eintr\"age haben.
  \emph{Entscheidungsm\"oglichkeiten} werden in der Sprache der
  Mathematik durch solche Variablen repr\"asentiert.
\item Wie ermittelt man die \emph{Kosten} einer TSP-Tour aus ihrer
  Adjazenzmatrix? Dies f\"uhrt auf die Matrix aller paarweisen
  Entfernungen, die manchen aus Atlanten bekannt sein d\"urfte.  Das
  \emph{Ziel}, eine m\"oglichst kurze Tour zu finden, wird repr\"asentiert
  durch eine \emph{Zielfunktion} in den Variablen.
\item Wie findet man heraus, ob eine gegebene Adjazenzmatrix wirklich
  alle \emph{Bedingungen} f\"ur eine TSP-Tour erf\"ullt? Dies ist der
  logisch anspruchsvollste Schritt und ein Sieg der Mathematik: Man
  kann dies f\"ur eine allgemeine TSP-Tour ausdr\"ucken, beschrieben
  durch mathematisch formulierte \emph{Restriktionen} f\"ur die
  Variablenbelegung (\emph{Gleichungen}, \emph{Ungleichungen},
  \emph{Ganzzahligkeitsbedingungen} in den Variablen). In der
  Mathematik kann man mit etwas rechnen, das man noch gar nicht kennt!
\end{enumerate}

Ergebnis dieser Schritte ist ein vollst\"andiges Modell f\"ur das TSP aus
der Klasse der \emph{Ganzzahligen Linearen Optimierungsaufgaben
  (ILP)}.  F\"ur diese gibt es Software, in der (fast) das ganze
mathematische Wissen \"uber solche Aufgaben kondensiert vorliegt.  Diese
Software wei\ss  weder etwas vom Handlungsreisenden noch vom ADAC.  Sie
kennt nur die Sprache der Mathematik.  Und f\"ur das Modell findet sie
L\"osungen f\"ur alle nicht zu gro\ss en Problembeispiele.  Wie gro\ss  {\glqq}nicht
zu gro\ss{\grqq} sein kann, zeigen die aktuellen Weltrekorde mit beweisbar
optimalen Touren durch \"uber zwanzigtausend
St\"adte.\footnote{\,\url{http://www.tsp.gatech.edu/}}

Im ADAC-Problem f\"uhrte dieses Vorgehen in der Praxis zum Erfolg.  Aber
die {\glqq}Rahmenhandlung{\grqq} zum ADAC-Problem wird auf einem narrativen Niveau
zum Abschluss gebracht.  Nat\"urlich nicht ohne einen Hinweis darauf,
dass man die Mathematik, die in der Optimierungs-Software steckt, an
der Universit\"at auch lernen kann (im Optimierungslabor jedoch leider
nicht).

Dass die \emph{Modellierung} mit ein wenig Hilfestellung gar nicht so
schwer ist, erleben die Sch\"ulerinnen und Sch\"uler beim Bearbeiten der
nun folgenden Modellierungsprojekte. Unser Vorgehen beschreiben
wir im Folgenden stellvertretend an einem R\"atsel aus dem Schach.

\section{Ein Modellierungsprojekt: nicht-schlagende Damenkonfigurationen im Schach}
\label{sec:Schachdamen}
Kann man acht Damen so auf einem Schachbrett platzieren, dass keine eine andere
schlagen kann, wenn sie sich, wie \"ublich, beliebig weit horizontal,
vertikal und diagonal \"uber das Brett bewegen d\"urfen?

Um dieses Schachr\"atsel als ILP modellieren zu k\"onnen, fragen wir nach den wesentlichen
Aspekten des Problems: Was kann man entscheiden? Nat\"urlich die Positionen der 8 Damen.
Der n\"achste Schritt ist die \"Ubersetzung unserer Entscheidungen in die Sprache der
Mathematik: Wie kann man diese \textit{kodieren}?

\begin{figure}[htp]
\begin{minipage}[t]{.43\textwidth}
    \centering
    \includegraphics[height=7\baselineskip]{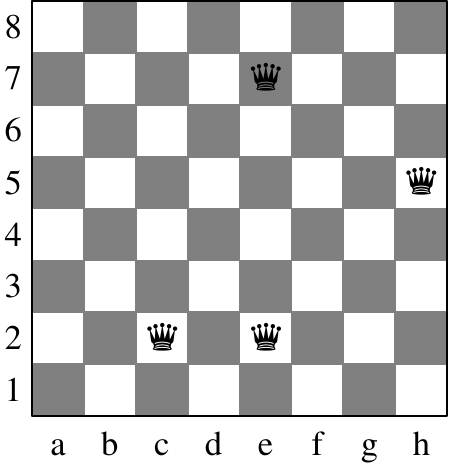}
    \quad
    \includegraphics[height=7\baselineskip]{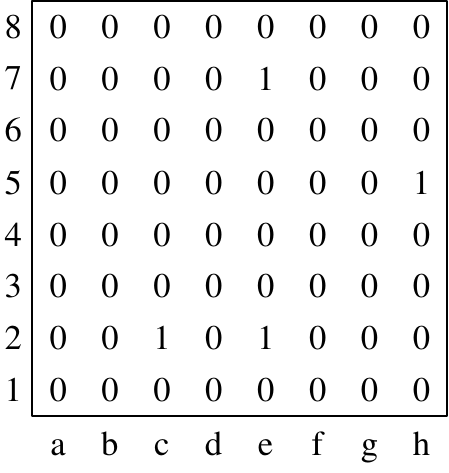}
    \caption{Gewohnte und \"ubersetzte Darstellung.}\label{fig:queens_rep}
\end{minipage}
\hspace*{.0\textwidth}
\begin{minipage}[t]{.54\textwidth}
    \centering
    \includegraphics[height=7\baselineskip]{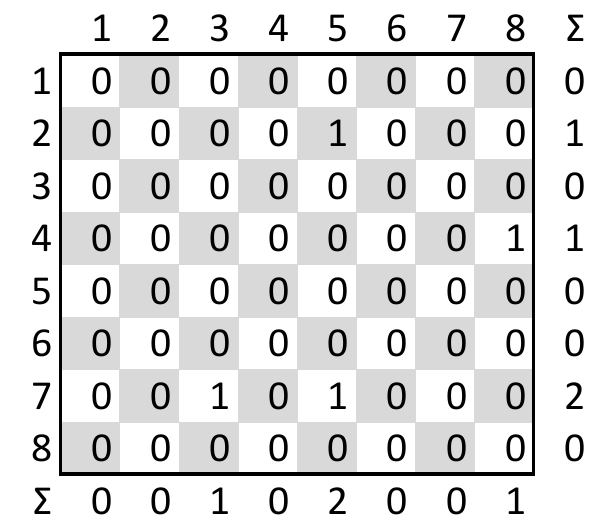}
    \quad
    \includegraphics[height=7\baselineskip]{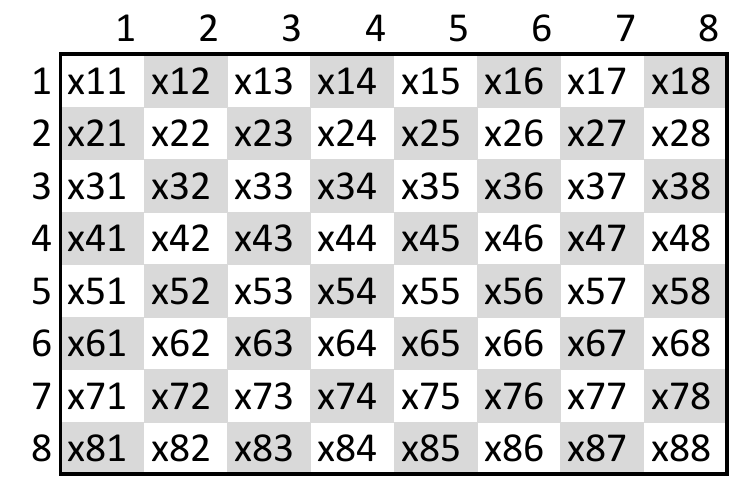}
    \caption{In einer Tabellenkalkulation.}\label{fig:queens_excel}
\end{minipage}
\end{figure}

Abbildung~\ref{fig:queens_rep} zeigt links eine Stellung von vier
Damen in gewohnter Manier. Rechts daneben haben wir auf schm\"uckendes Beiwerk
verzichtet: Felder mit einer Dame werden durch eine $1$, Felder ohne
durch eine $0$ ausgedr\"uckt. Diese Darstellung hat den Vorteil,
dass wir damit -- z.\,B.{ }in einer Tabellenkalkulation wie Excel -- rechnen
k\"onnen.

Wo die Damen letztlich stehen m\"ussen, wissen wir aber noch gar nicht, lediglich
\emph{dass} auf jedem Feld eine steht oder nicht. Das dr\"ucken wir durch Variablen aus,
die nur die Werte $0$ oder $1$ annehmen k\"onnen und deren Name Auskunft \"uber
die Position des zugeh\"origen Feldes gibt,
siehe Abbildung~\ref{fig:queens_excel} rechts. Der dabei vorgenommene
Wechsel von der gewohnten
alphanumerischen Feldbezeichnung zu einer rein numerischen wird sich sp\"ater
als n\"utzlich herausstellen.

Es k\"onnte sein, dass es nicht m\"oglich ist, acht Damen
(nicht-schlagend) auf einem Schachbrett zu platzieren.
Vorsichtshalber formulieren wir die Aufgabe deshalb um zu
{\glqq}Platziere m\"oglichst viele Damen, die sich gegenseitig nicht schlagen
k\"onnen{\grqq}. Deren Anzahl l\"asst sich dann bestimmen, indem man
die $0/1$-Werte aller 64~Felder aufaddiert. Unser \emph{Ziel} muss es damit sein, diese
Summe m\"oglichst gro\ss  zu machen. Das Sch\"one an Mathematik:
Diese Idee funktioniert auch, wenn wir Variablen
(also noch unbekannte Gr\"o\ss en) verwenden.
Unsere \textit{Zielfunktion} lautet dann
\begin{equation*}
  \max \quad x_{1,1} + x_{1,2} + x_{1,3} + x_{1,4} + x_{1,5} +
             x_{1,6} + x_{1,7} + x_{1,8} + x_{2,1} + \dots + x_{8,7} + x_{8,8}.
\end{equation*}

Nun m\"ussen wir daf\"ur sorgen, dass die {\glqq}Regeln{\grqq} eingehalten werden, sich die
platzierten Damen also nicht gegenseitig schlagen k\"onnen.
In Abbildung~\ref{fig:queens_rep} ist das verletzt, leicht zu erkennen an den
Zweiern in Abbildung~\ref{fig:queens_excel} links, welche die Zeilen- und
Spaltensummen angeben. Eine Stellung ist nur dann zul\"assig, wenn diese Summe
$0$ oder $1$ ist. Ausgedr\"uckt mit unseren \emph{Entscheidungsvariablen} lautet
diese \emph{Nebenbedingung} f\"ur Zeile~7:
\begin{align*}
x_{7,1} + x_{7,2} + x_{7,3} + x_{7,4} + x_{7,5} + x_{7,6} + x_{7,7} + x_{7,8} \leq 1.
\end{align*}

\begin{figure}[htp]
  \begin{center}
    \includegraphics[height=7\baselineskip]{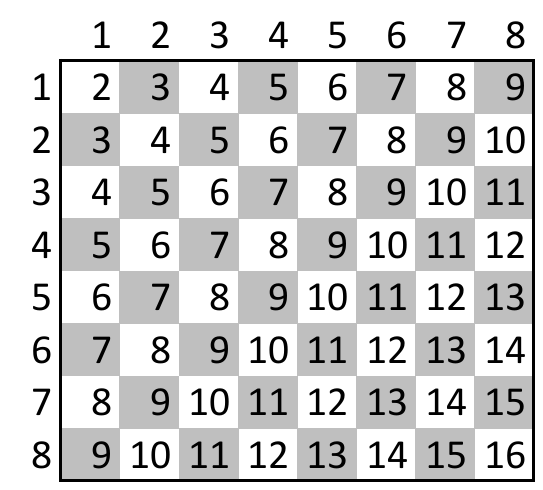}
    \quad
    \includegraphics[height=7\baselineskip]{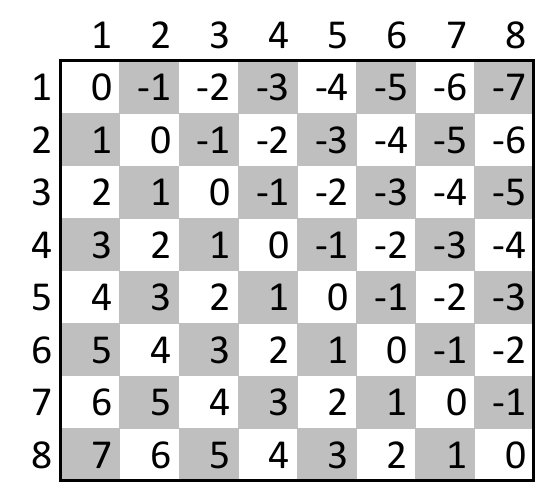}
    \quad
    \includegraphics[height=7\baselineskip]{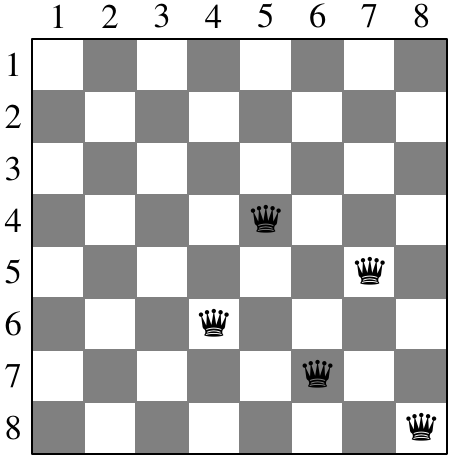}
    \quad
    \includegraphics[height=7\baselineskip]{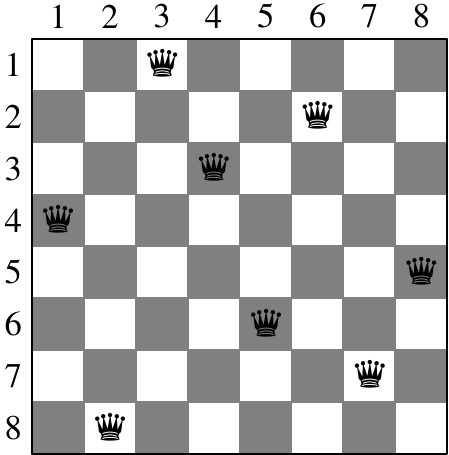}
  \end{center}
  \caption{Summe bzw.~Differenz der Feldindizes und zwei Beispielkonfigurationen.}
  \label{fig:damen_index}
\end{figure}

Damit bleiben die Diagonalen. Diese lassen sich ganz analog behandeln. F\"ur die
Hauptdiagonale von links oben nach rechts unten lautet unsere Bedingung
\begin{equation*}
    \fontsize{7}{7}\selectfont
x_{\underset{1-1=0}{\underbrace{1,1}}} + x_{\underset{2-2=0}{\underbrace{2,2}}}
+ x_{\underset{3-3=0}{\underbrace{3,3}}} + x_{\underset{4-4=0}{\underbrace{4,4}}}
+ x_{\underset{5-5=0}{\underbrace{5,5}}} + x_{\underset{6-6=0}{\underbrace{6,6}}}
+ x_{\underset{7-7=0}{\underbrace{7,7}}} + x_{\underset{8-8=0}{\underbrace{8,8}}}
\leq 1.
\end{equation*}
Die geschweiften Klammern verdeutlichen, dass man sich -- dank der numerischen
Feldbezeichner -- eines Kniffs bedienen kann: Bei Aufw\"artsdiagonalen ist jeweils
die \emph{Summe} der Indizes gleich, bei Abw\"artsdiagonalen die \emph{Differenz} der Indizes,
siehe Abbildung~\ref{fig:damen_index}. Damit lassen sich die Felder einer
Diagonale kompakt beschreiben.

Insgesamt ergeben
sich je acht Nebenbedingungen f\"ur Zeilen bzw.{ }Spalten und $2\cdot15$
Nebenbedingungen f\"ur die Diagonalen. Mit etwas mathematischem Formalismus
und durch Ausnutzen des eben gesehenen Musters kann man sich einige Schreibarbeit
sparen und das ILP-Modell kompakt hinschreiben:

\begin{align*}
 \max \sum\nolimits_{i=1}^8 \sum\nolimits_{j=1}^8 x_{i,j} \qquad &\textrm{unter den
 Nebenbedingungen}\\
 x_{i,j}\in\{0,1\}            \qquad &\forall i, j = 1, \dots, 8,\\
 \sum\nolimits_{j=1}^8 x_{i,j}  \leq 1 \qquad &\forall i = 1, \dots, 8,\\
 \sum\nolimits_{\,i=1}^8 x_{i,j}  \leq 1 \qquad &\forall j = 1, \dots, 8,\\
 \sum\nolimits_{i+j=k} x_{i,j}  \leq 1 \qquad &\forall k = 2, \dots, 16,\\
 \sum\nolimits_{i-j=k} x_{i,j}  \leq 1 \qquad &\forall k = -7, \dots, 7.
\end{align*}

\"Ubergibt man dieses Modell an ein Computerprogramm zum L\"osen von
ILPs\footnote{\,z.\,B.{ }die ZIB Optimization Suite, siehe \url{http://zibopt.zib.de/}.},
k\"onnte man die rechte Konfiguration aus
Abbildung~\ref{fig:damen_index} erhalten (insgesamt gibt es 92
L\"osungen).  Die Konfiguration links daneben mit nur f\"unf Damen ist auch
etwas Besonderes. Warum?\footnote{\,Die Aufl\"osung findet sich am Ende des Artikels.}

Selten gibt es in der Mathematik nur einen Weg zum Ziel, und so lassen sich auch
beim Damenproblem alternative Modelle finden, siehe dazu~\cite[S.\,32\,ff.]{koch:queens}.

\section{Weitere Modellierungsprojekte}
\label{sec:AndereModelle}
Den Sch\"ulerinnen und Sch\"ulern stellen wir viele weitere Aufgabenstellungen zur
Wahl. Das Projekt der nicht-schlagenden Damen haben wir hier genauer beleuchtet,
da es die \"Ubersetzung von Entscheidungsm\"oglichkeiten, Zielen und
Restriktionen eines direkt nachvollziehbaren Problems in Variablen,
Bewertungsfunktionen, Gleichungen und Ungleichungen veranschaulicht.

Im Laufe der Zeit haben sich einige Modellierungsprojekte herauskristallisiert,
deren Problemstellungen wir im Folgenden beschreiben wollen. Auswahlkriterien
waren die potentielle \textit{N\"ahe} zum Alltag der Sch\"uler (abh\"angig von
der Interessenslage), die M\"oglichkeit, ohne viel theoretische Vorkenntnisse
aktiv mitmachen zu k\"onnen, eine \"uberschaubare Komplexit\"at (bedingt durch
zeitliche Einschr\"ankungen und den anvisierten motivierenden Charakter) und
teilweise gr\"o\ss ere konzeptionelle \"Uberschneidungen zwischen den einzelnen
Projektvorschl\"agen, um Erfolgserlebnisse durch \textit{Transferleistungen} zu
erm\"oglichen. Des Weiteren sollte das Problem nicht allzu einfach \textit{mit
der Hand} zu l\"osen sein. Die nachfolgend vorgestellte Auswahl ist keineswegs
abschlie\ss end zu sehen. Sie soll dem Leser als Anregung dienen und die Breite
m\"oglicher Problemstellungen andeuten. Ein letzter Warn-
bzw.{ }Motivationshinweis: Die Grenze zwischen reinen Spielproblemen, \"okonomisch
relevanten Problemen und Forschungsproblemen ist flie\ss end.

\begin{figure}[htp]
\begin{minipage}[t]{.49\textwidth}
    \centering
    \includegraphics[height=7\baselineskip]{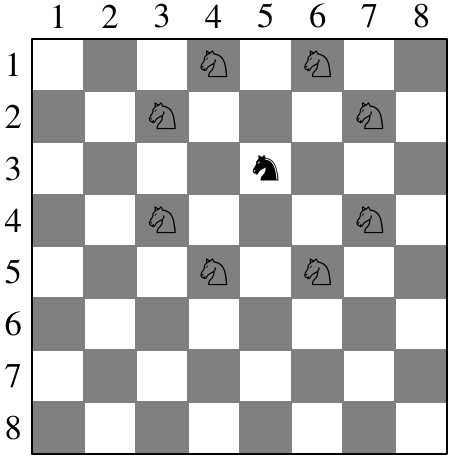}
    \quad
    \includegraphics[height=7\baselineskip]{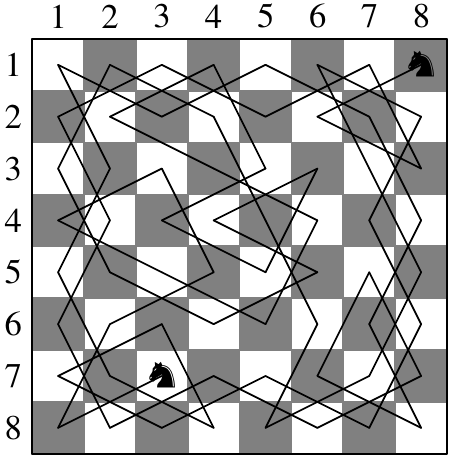}
    \caption{Links: Der R\"osselsprung. Rechts: Eine offene Springertour.}%
    \label{fig:knight_moves}
\end{minipage}
\hspace*{.0\textwidth}
\begin{minipage}[t]{.49\textwidth}
    \centering
    \includegraphics[height=7\baselineskip]{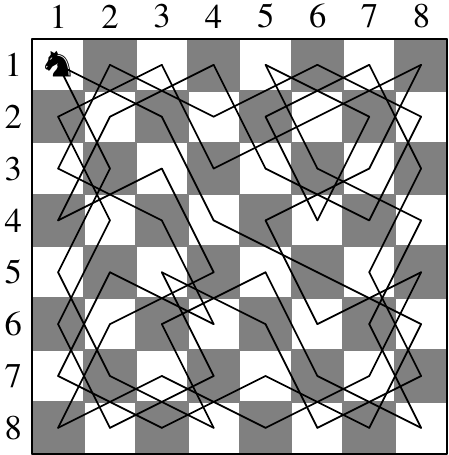}
    \quad
    \includegraphics[height=7\baselineskip]{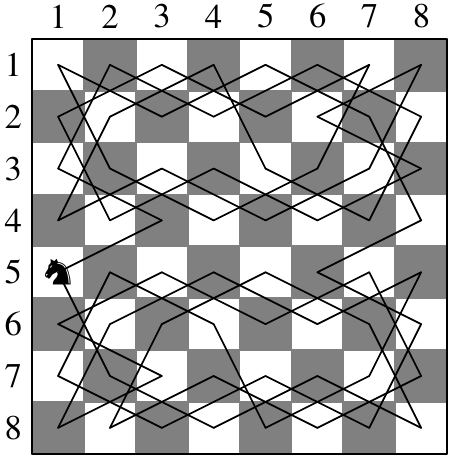}
    \caption{Zwei verschiedene geschlossene Springertouren.}
\end{minipage}
\end{figure}

Bleiben wir zun\"achst beim Schach. Eine \textit{Springertour} ist eine
Route auf einem leeren Schachbrett, bei der jedes Feld genau ein Mal besucht
wird. Bewegen darf sich ein Springer jeweils zwei Felder gerade aus und eines
zur Seite, siehe Abbildung~\ref{fig:knight_moves}. Sind Start- und Endfeld einen
Springerzug voneinander entfernt, so spricht man von einer
\textit{geschlossenen} Springertour. Beginnend mit dem Schweizer Mathematiker
Leonard Euler haben sich seit 1759 viele Mathematiker und Hobby-T\"uftler mit
dem Problem besch\"aftigt. So wurde die Frage nach der algorithmischen
Konstruktion solcher Touren auf verallgemeinerten $n\times n$-Schachbrettern
auch im Rahmen einer "`Jugend forscht"'-Arbeit gekl\"art,
siehe~\cite{0793.68113} f\"ur Details.

Eine etwas \textit{praktischere} Problemstellung ist \textit{Routenplanung} --
Navigationsger\"ate f\"ur den Stra\ss enverkehr finden sich heutzutage in fast jedem
PKW. Nach Eingabe des Start- und Zielorts wird auf Knopfdruck eine k\"urzeste
oder schnellste
Route berechnet. Kennt man die Distanz bzw.{ }den Zeitverbrauch zwischen
allen m\"oglichen benachbarten Zwischenzielen, l\"asst
sich auch diese Aufgabe mit Hilfe ganzzahliger linearer Optimierung
modellieren.\footnote{\,Navigationsger\"ate verwenden allerdings ma\ss
geschneiderte K\"urzeste-Wege-Algorithmen f\"ur dieses Problem, die deutlich
schneller sind. Stichworte sind: Dijkstra oder A$^\star$-Algorithmus.}

\begin{figure}[htp]
\begin{minipage}[t]{.4\textwidth}
    \centering
    \includegraphics[height=7\baselineskip]{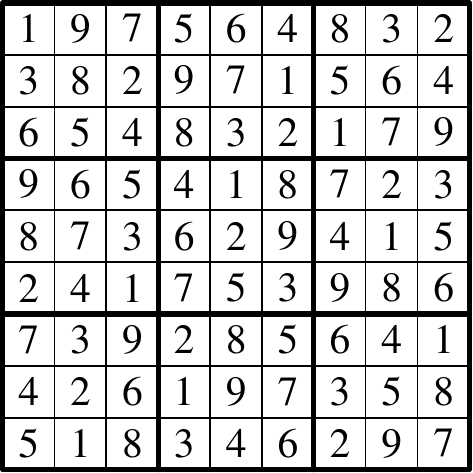}
    \caption{Gel\"ostes Sudoku.}\label{fig:sudoku_solved}
\end{minipage}
\hspace*{.0\textwidth}
\begin{minipage}[t]{.55\textwidth}
    \centering
    \includegraphics[height=7\baselineskip]{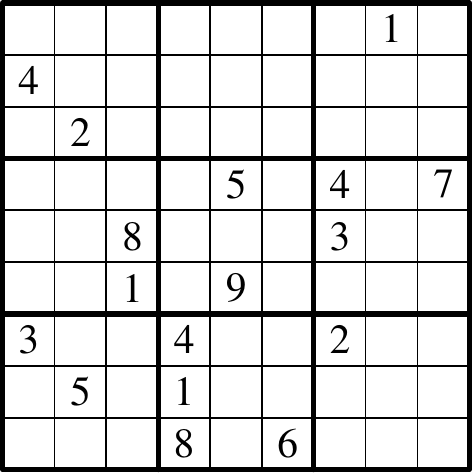}
    \quad
    \includegraphics[height=7\baselineskip]{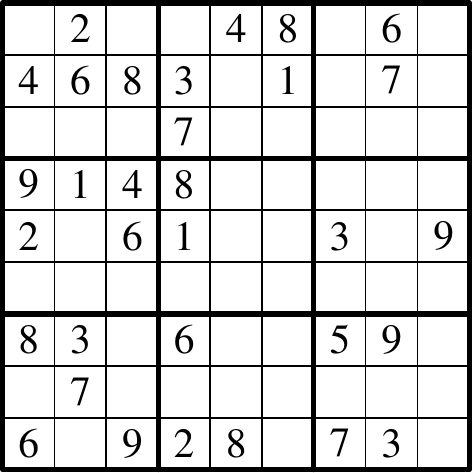}
    \caption{Ausgangssituationen.}\label{fig:sudoku_unsolved}
\end{minipage}
\end{figure}

Mit \emph{Sudoku} haben wir einen weiteren Leckerbissen f\"ur R\"atselknacker im
Programm.  Das Zahlenpuzzle ben\"otigte seit seiner Erfindung 1979\footnote{\,Howard
Garns erfand es seinerzeit unter dem Namen \emph{NumberPlace}, siehe
\url{http://de.wikipedia.org/wiki/Sudoku}.} einige Jahre und den Umweg \"uber
Japan (daher der Name) bis es hierzulande gro\ss e Beliebtheit erlangte. Dabei soll
ein $9\times 9$ Gitter so mit den Zahlen von $1$ bis $9$ ausgef\"ullt werden, dass
in jeder Spalte, jeder Zeile und in jedem Block ($3\times 3$-Untergitter) jede
Zahl genau ein Mal vorkommt, wobei stets einige Zahlen vorgegeben sind, siehe
Abbildung~\ref{fig:sudoku_unsolved}. Die Suche nach einer \textit{zul\"assigen}
Sudokul\"osung kann wiederum als ILP modelliert werden.\footnote{\,Das Ergebnis
(Abbildung~\ref{fig:sudoku_solved}) ist ein {\glqq}alter Hut{\grqq}, ein
lateinisches Quadrat -- auch damit befasste sich schon Leonard Euler. Als ILP
wird es in~\cite{koch:sudoku} und~\cite{bartlett+etal:sudoku} betrachtet.} Noch etwas
interessanter wird es, wenn man zus\"atzlich fordert, dass es nur eine, und damit
eindeutige, L\"osung geben darf.\footnote{\,Als der
Sudoku-Hype noch gr\"o\ss er war, gab es einige wohl sehr eilig produzierte
R\"atselhefte mit nicht eindeutig l\"osbaren Sudokus -- einer der Autoren hatte
sich damals beim zust\"andigen Verlag beschwert und zumindest ein weiteres
Gratisheft erhalten. Eines der beiden unvollst\"andig ausgef\"ullten Sudokus aus
Abbildung~\ref{fig:sudoku_unsolved} ist nicht eindeutig l\"osbar.} Die Frage
nach der minimal n\"otigen Anzahl an Hinweisen, bei der ein Sudoku eine eindeutige
L\"osung hat, konnte erst vor Kurzem beantwortet werden: 17~ausgef\"ullte Felder
werden mindestens ben\"otigt.\footnote{Ein Forscherteam um Gary McGuire konnte
durch vollst\"andige Enumeration zeigen, dass kein eindeutig l\"osbares Sudoku
mit nur 16~Hinweisen existiert, siehe \url{http://www.arxiv.org/abs/1201.0749}
bzw.{} \url{http://www.math.ie/checker.html}.}

Ein weiteres, recht anschauliches Problem, welches in vielen Modellierungen als
Teilaspekt vorkommt, ist das sogenannte \textit{Rucksackproblem}. Hierbei hat man
eine Menge von Gegenst\"anden, beispielsweise ein Handy, einen MP3-Player,
Schokolade, ein Tagebuch, u.\,v.\,a.\,m.{ }zur Auswahl, die alle ein gewisses Gewicht und einen
pers\"onlich festgelegten Nutzen haben. Unter einer Gewichtsrestriktion m\"ochte
man nun den Rucksack so best\"ucken, dass die summierten Nutzenwerte maximal
sind.

\begin{figure}[htp]
    \centering
    \includegraphics[height=4\baselineskip]{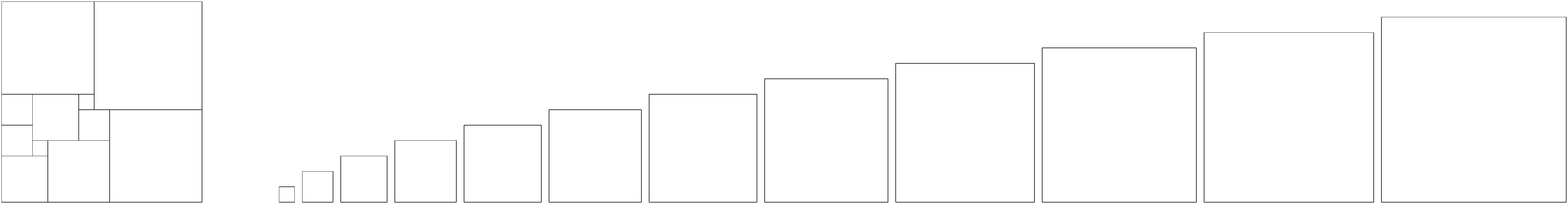}
    \caption{Fliesplan f\"ur ein quadratisches $13\times 13$ Zimmer und zur
    Auswahl stehende Fliesen.}\label{fig:tiling}
\end{figure}

Zum Abschluss\footnote{\,Weitere Problemstellungen wie
  \textit{Globetrotter} oder \textit{Blind-Dance} befinden sich unter
  \url{http://www.wm.uni-bayreuth.de/index.php?id=optlabor} bzw.{
  }\cite{Kombinatorische_Optimierung_erleben}.} m\"ochten wir ein
geometrisches Fliesenproblem erw\"ahnen. Ein $13\times 13$-Zimmer soll
mit kleineren $a\times a$-Fliesen ($a = 1, 2, 3, \dots, 12$) vollst\"andig (und
\"uberlappungsfrei) gefliest werden. Eine M\"oglichkeit mit 11~Fliesen
ist in Abbildung~\ref{fig:tiling} dargestellt. Diese L\"osung finden
und zeigen, dass es keine mit weniger Fliesen gibt, kann man mit Hilfe
einer ILP-Modellierung.

\section{Unser Weg zum Material}
\label{sec:Material}
Unser Optimierungslabor war urspr\"unglich gedacht als ein Miniworkshop
zur computergest\"utzten L\"osung von (vereinfachten) Alltagsproblemen.
Blaupause war eine Blockvorlesung f\"ur Studierende im Umfang einer
vierst\"undigen Vorlesung mit einer zweist\"undigen \"Ubung: eine
Tour-de-force, aber bei den Studierenden erfolgreich.

Nun waren die Voraussetzungen in einer eint\"agigen Veranstaltung
\"uber einige Stunden anders. Trotzdem stand anfangs die
Computernutzung im Zentrum. Gedanken gemacht hat sich die
teilnehmende Gruppe dann vor dem Bildschirm in Zweier- oder
maximal Dreierteams. Wir haben dann festgestellt, dass dies f\"ur
die eigentlich intellektuell anspruchsvolle Aufgabe, die Logik
eines Problems zu erfassen und mathematisch zu formalisieren,
keine attraktive Arbeitsumgebung ist. Ferner stellte das
Selbereintippen des Modells f\"ur die Optimierungssoftware zwar eine
Aktivierung der Gruppe dar, verlagerte aber den Schwerpunkt zu
sehr auf das rein Handwerkliche der Computernutzung. Die
unterschiedliche Affinit\"at zum Umgang mit dem Computer selbst
stellte ein weiteres Problem dar: Informatikbegeisterte fanden
sich schnell zurecht, andere waren dadurch eher \"uberfordert, was
der Motivation insgesamt abtr\"aglich war.

In einem Zwischenschritt haben wir Arbeitsbl\"atter entworfen, die durch
gezielte Aufgaben zum Selberl\"osen die Gruppe auch ohne l\"uckenlose
Betreuung zum Modell hinf\"uhren sollte.  Alles in allem blieb aber der
Betreuungsaufwand immens, und der Computer lud alle die zum
Abschweifen ein, bei denen gerade kein Betreuungsgespr\"ach lief.

Wir entschlossen uns dann zu einer klareren Schwerpunktsetzung und
einer deutlicheren Abgrenzung zur mehrt\"agigen Blockveranstaltung f\"ur
Studierende: Die \"Ubersetzung eines Problems in die Sprache der
Mathematik sollte durch Material und Sitzordnung in den Mittelpunkt
r\"ucken. Da die L\"osung durch die Software am Ende nicht verzichtbar
ist, haben wir diesen {\glqq}kr\"onenden{\grqq} Abschluss ins Plenum verlegt: Ein
erfahrenes Teammitglied bedient am Schluss den Computer, dessen
Bildschirmausgabe auf dem Beamerbild live verfolgt werden kann.
\begin{figure}[htp]
    \centering
    \begin{minipage}[t]{.227\textwidth}
        \includegraphics[height=24\baselineskip]{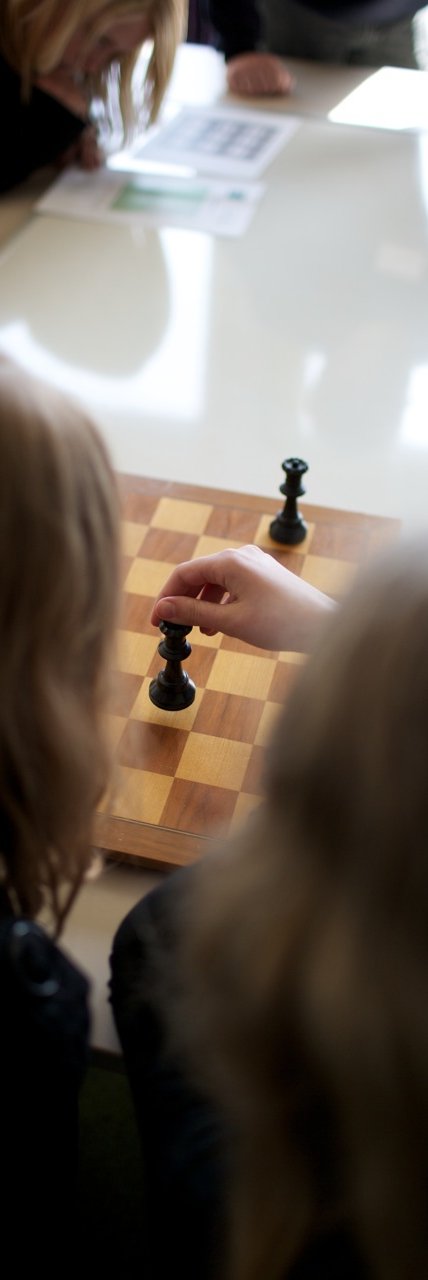}
    \end{minipage}
    \begin{minipage}[b]{.56\textwidth}
        \includegraphics[height=12\baselineskip]{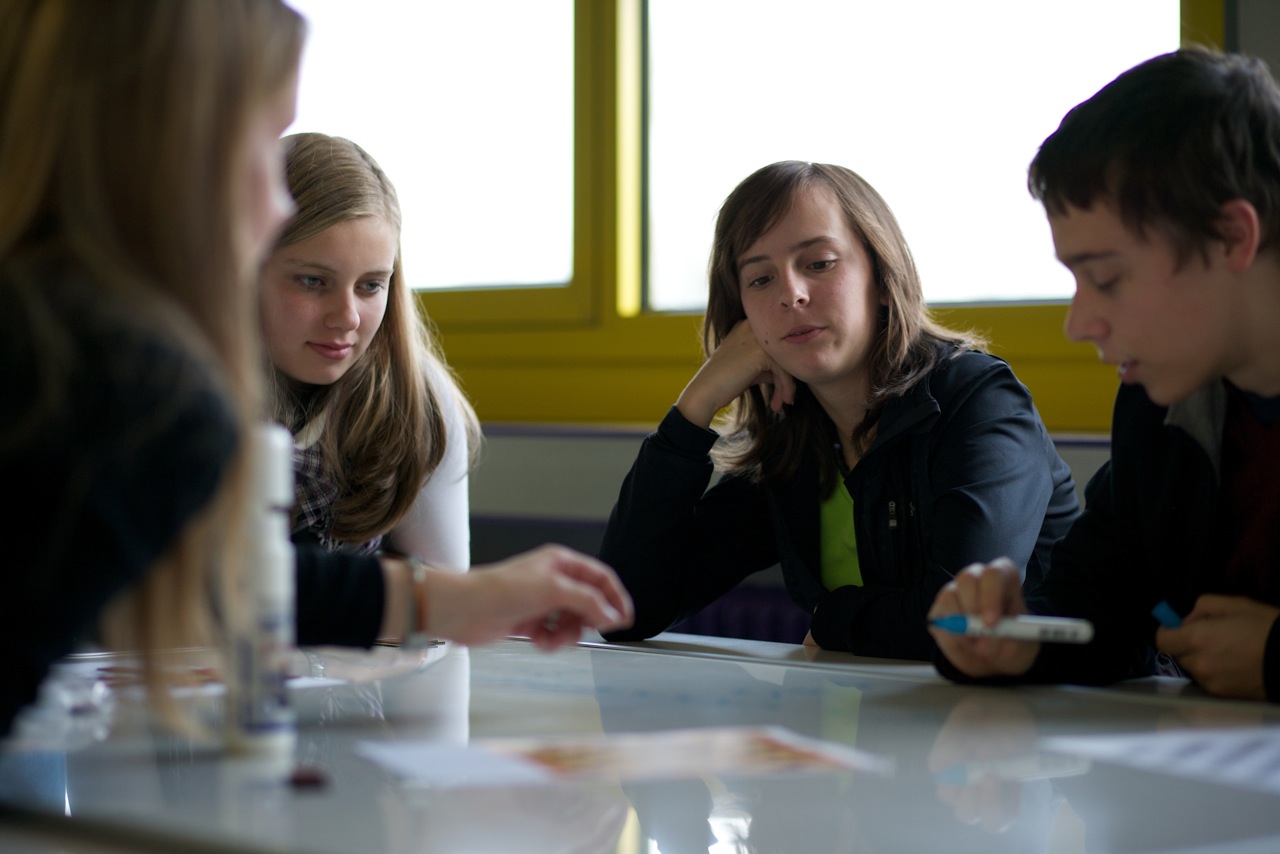}
        \includegraphics[height=18\baselineskip,angle=90]{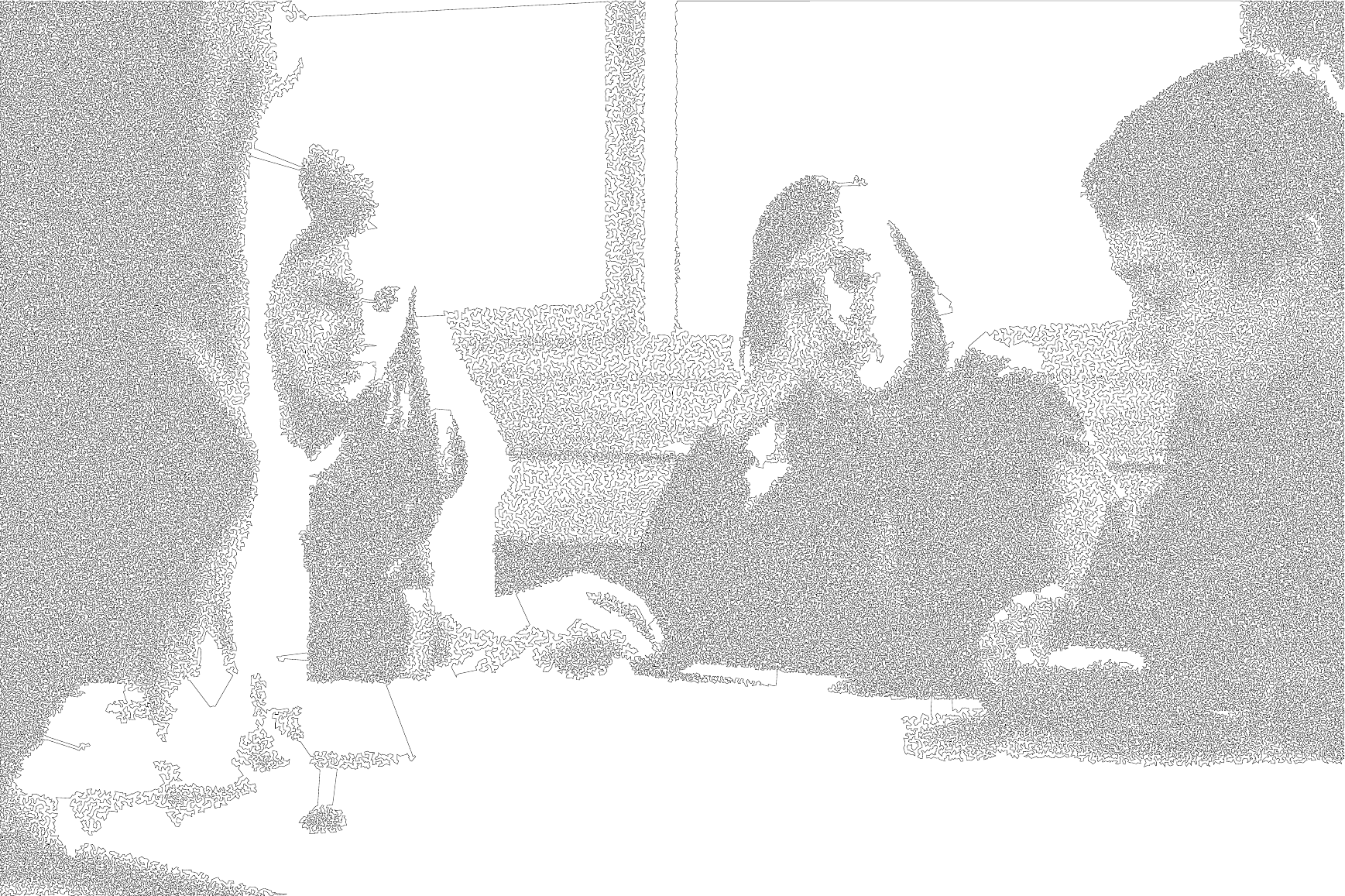}
    \end{minipage}
    \caption{Links: W\"ahrend der Probierphase. Daneben: Teamdiskussion als Foto
    und darunter als TSP-Tour mit 447\,256 zu besuchenden {\glqq}St\"adten{\grqq}; Software
    mit freundlicher Genehmigung des Zuse Institut in Berlin.}
    \label{fig:Diskussion_TSP}
\end{figure}
Das Erarbeiten der Modelle wird nun mit dreidimensionalem Material
unterst\"utzt: F\"ur jedes Team wird ein gro\ss er Tisch mit einer
Magnettafel bedeckt.  Auf dieser Tafel werden reale Objekte der
Projekte (z.\,B.{ }ein Schachbrett mit acht Damen, eine Waage etc.),
beschreibbare Magnete, Whiteboardmarker, Papier etc. bereitgestellt.
Von allen Seiten kann am Objekt probiert werden, Magnete beschriftet,
platziert und verschoben werden, auf die Tafel geschrieben und
gewischt werden.  F\"ur eine Sechsergruppe und eine Betreuungsperson ist
problemlos simultane Aktivit\"at und Diskussion m\"oglich.

\begin{figure}[htp]
    \centering
    \begin{minipage}[b]{.426\textwidth}
        \includegraphics[height=10\baselineskip]{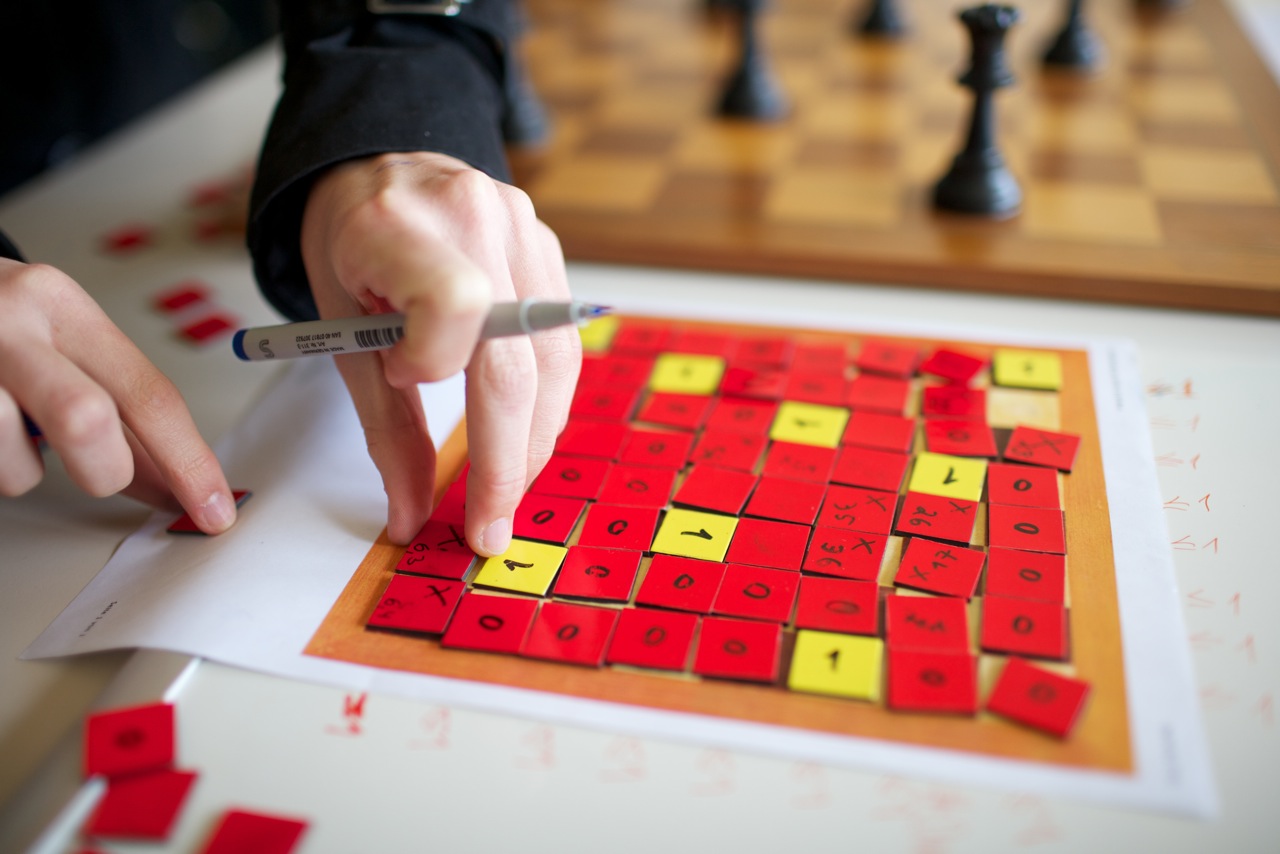}\vfill
        \includegraphics[height=10\baselineskip]{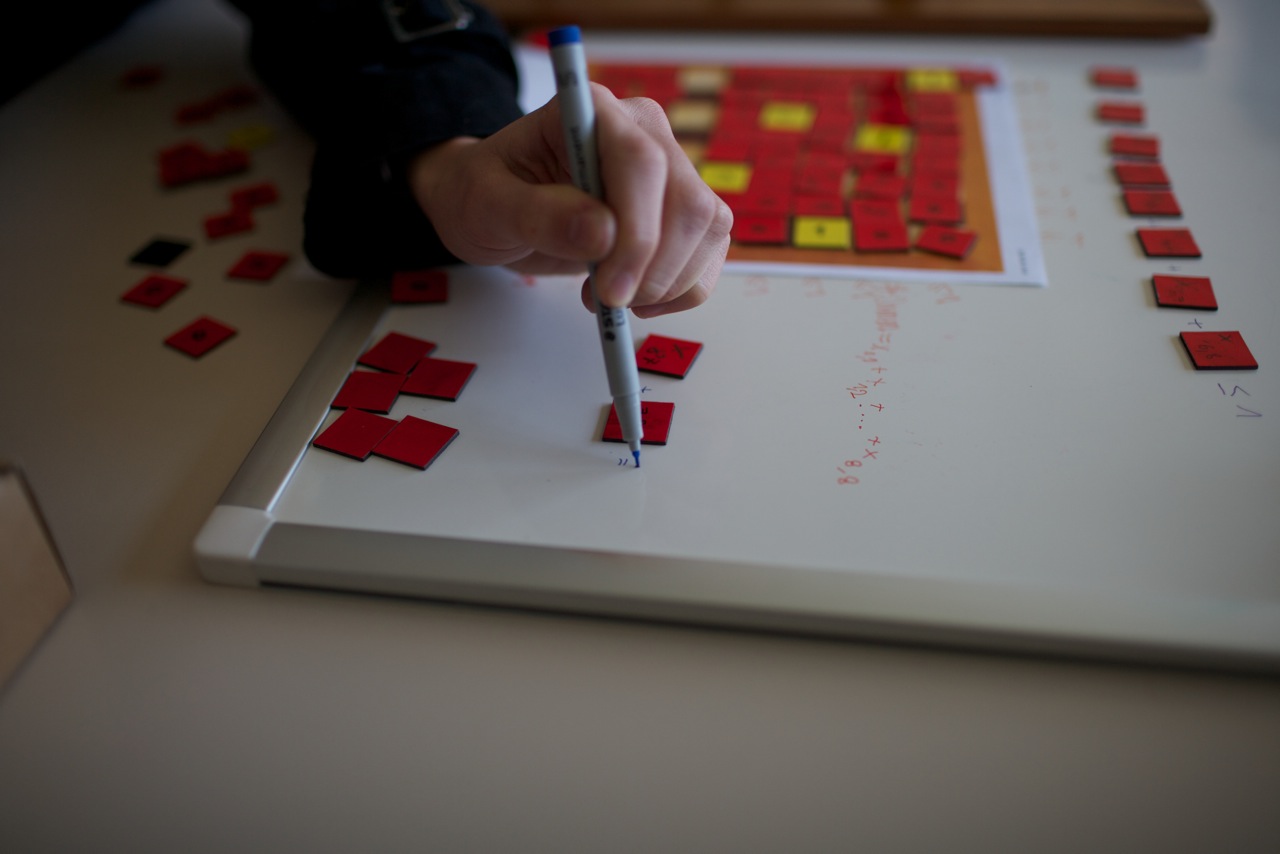}
    \end{minipage}
    \begin{minipage}[t]{.4\textwidth}
        \includegraphics[height=20\baselineskip]{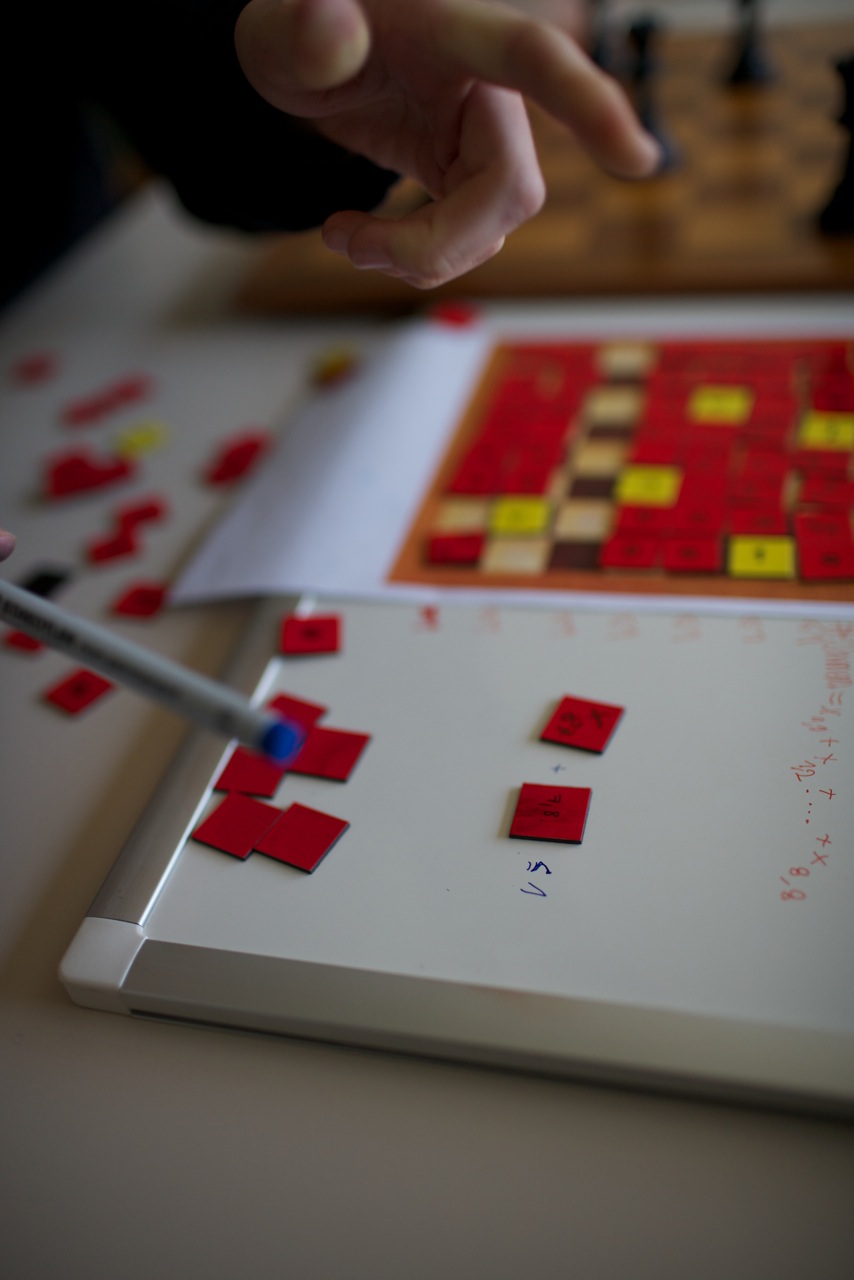}
    \end{minipage}
    \caption{Von 0/1-Belegungen \"uber Variablen hin zu Ungleichungen.}\label{fig:magnets}
\end{figure}

Mit den beschreibbaren Magneten k\"onnen L\"osungsvarianten (z.\,B.{ }eine
Damenkonfiguration auf dem Schachbrett) kodiert werden (siehe
Abbildung~\ref{fig:magnets}), indem man die
Einsermagnete auf die Felder mit den Damen legt und die Nullermagnete
auf die restlichen Felder. Variablen k\"onnen zun\"achst den
Entscheidungen r\"aumlich zugeordnet werden (z.\,B.{ }durch Auflegen von
$x_{23}$ auf das Schachfeld C7).  Danach k\"onnen die identischen
Magnete zur Synthese einer Formel auf die Tafel verfrachtet werden.
Ergebnissicherung kann dann klassisch durch Abschreiben erfolgreicher
L\"osungsans\"atze von Tafel auf Papier erfolgen. Die r\"aumliche N\"ahe aller
beteiligten Personen und Objekte sowie die gro\ss e Aktionsfl\"ache hatten
einen enorm positiven Einfluss auf die Mitwirkung.

\begin{figure}[htp]
\centering
  \includegraphics[height=7\baselineskip]{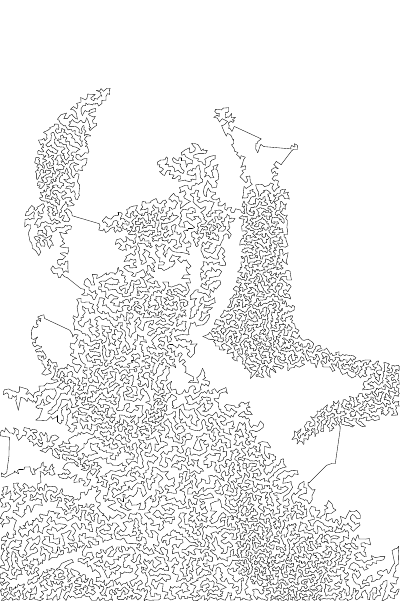}
  \includegraphics[height=7\baselineskip]{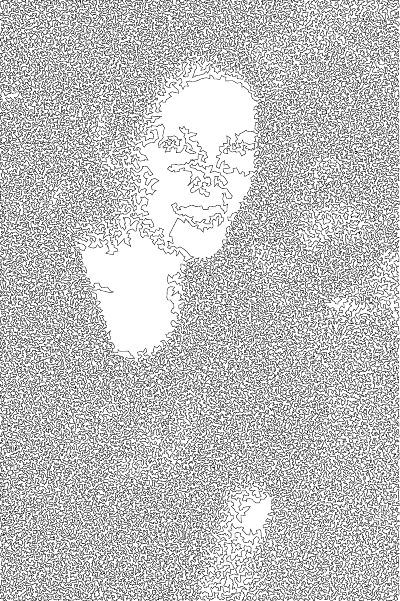}
  \includegraphics[height=7\baselineskip]{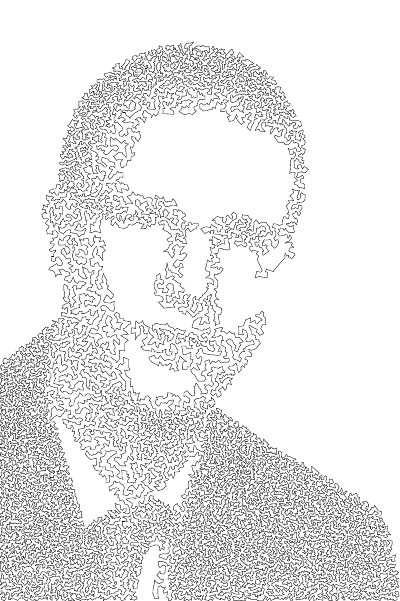}
  \includegraphics[height=7\baselineskip]{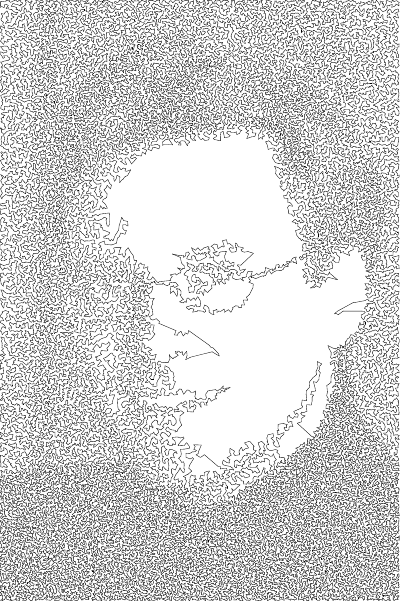}
  \includegraphics[height=7\baselineskip]{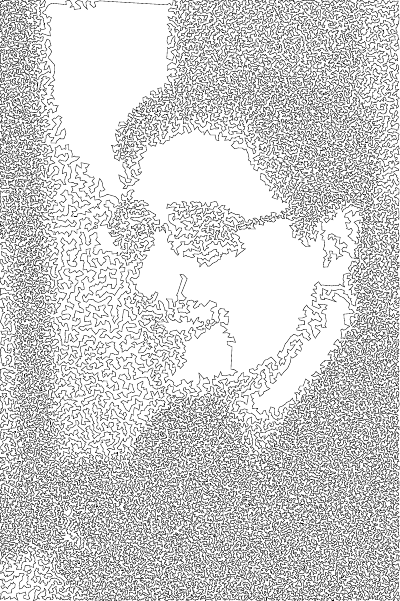}
  \includegraphics[height=7\baselineskip]{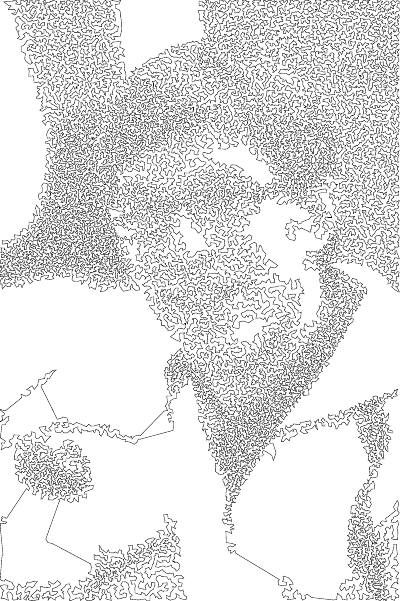}
  \caption{Eine TSP-Tour als Erinnerung exemplarisch am \glqq Optimierungsteam \grqq,
        Software mit freundlicher Genehmigung des Zuse Institut in Berlin.}
  \label{fig:Team_TSP}
\end{figure}

Ein {\glqq}Mitgebsel{\grqq} zum Thema (z.\,B.{ }ein Portrait der Teilnehmer, in
Punkte aufgel\"ost und durch eine TSP-Tour verbunden, siehe
Abbildung~\ref{fig:Team_TSP}) hat eine Verst\"arkung des
Erinnerungseffekts zum Ziel.

\section{Gesammelte Erfahrungen}
Die eben beschriebene Form des Optimierungslabors haben wir zu zahlreichen
Gelegenheiten angeboten und dabei einige Erfahrungen gesammelt, \"uber die wir
hier berichten wollen.

Um aktiv an der Gruppenarbeit teilnehmen zu k\"onnen, sollten die Teilnehmenden
mit Grundlagen aus der linearen Algebra, insbesondere der L\"osung linearer
Gleichungssysteme bzw.{ }dem Prinzip der Darstellung von Zusammenh\"angen mittels
der Verwendung von Variablen vertraut sein. Als Zielgruppe f\"ur unser Angebot
sehen wir daher die neunte Klassenstufe, u.\,U.{ }die achte.


Gerade der Umschwung auf die Gruppenarbeit mit einer Betreuungsperson f\"ur sechs bis
acht Personen erwies sich als gro\ss e Verbesserung. So haben auch
mathematisch weniger versierte Sch\"uler(innen) die M\"oglichkeit, sogar nach einer
kurzen Durststrecke, zu einem sp\"ateren Zeitpunkt wieder aktiv in die
Zusammenarbeit einzusteigen. Als Beispiel seien hier Sudoku und die
Nebenbedingungen {\glqq}eine Vier pro Spalte{\grqq} und {\glqq}eine Vier pro Zeile{\grqq} genannt. Diese
sind insofern analog, als dass bei ersterer der Zeilenindex und bei letzterer der
Spaltenindex variabel ist. Sobald eine dieser {\glqq}\"ahnlichen{\grqq}
Nebenbedingungen entdeckt wurde, fand sich fast immer ein anderes
Gruppenmitglied, das die zweite Bedingung formulieren konnte.

Im Allgemeinen fiel es den Sch\"uler(inne)n leichter, Ideen zu \"au\ss ern
bzw.{ }Zielfunktion oder Nebenbedingungen m\"undlich zu formulieren,
als schriftlich auf der Tafel niederzuschreiben. Der spielerische
Einstieg durch Probieren (z.\,B.{ } L\"osen von Sudokus, Aufstellen der
acht Damen auf ein gro\ss es Schachbrett usw.) fand stets sehr gro\ss en
Anklang.\footnote{\,Diesen Aspekt des {\glqq}Vertrautwerdens{\grqq} mit dem Problem
  betonen wir, ist er doch unverzichtbarer Bestandteil des L\"osens an
  sich.} Das erkl\"arte Ziel die Modellprojekte so auszuw\"ahlen, dass sie
zum Ausprobieren einladen, aber nicht im Handumdrehen gel\"ost werden
k\"onnen, scheinen wir gut erreicht zu haben: So bekommt zum Beispiel
die Sudoku-Gruppe seit einiger Zeit zwei Sudokus, darunter ein
unl\"osbares. Taucht hier beim Selberl\"osen ein Fehler auf, ist es
schwer nachzuvollziehen, ob festgesetzte oder eingetragene Zahlen
verantwortlich sind -- f\"ur den Computer eine Sache von
Sekundenbruchteilen.\footnote{\,F\"ur viele Optimierungsprobleme gibt es
  speziell entwickelte L\"osungsalgorithmen, die schneller als ein
  geeigneter ILP-Ansatz sind.  Dies herauszustellen k\"onnte Gegenstand
  einer anderen Initiative sein, die gr\"o\ss tm\"ogliche Effizienz in den
  Mittelpunkt stellt.  Die Algorithmen m\"ussen dann allerdings auch
  etwas Spezielles \"uber das Problem wissen und brauchen die Daten in
  speziellerer Form.  Sie sind daher f\"ur andersartige Probleme nicht
  mehr zu gebrauchen.  Wir wollen in \emph{diesem} Projekt die
  Universalit\"at mathematischer Sprache illustrieren und zielen daher
  auf die logisch korrekte Anwendung einer mathematisch m\"oglichst
  universellen Methode.} Die abschlie\ss ende Demonstration der
\"Ubertragung des Modells auf den Computer und die Pr\"asentation der
Ergebnisse wird mit gro\ss em Interesse verfolgt. Insbesondere stellen
wir selten Schwierigkeiten beim R\"uck\"ubersetzen der Ergebnisse in den
Anwendungskontext fest.

Uns als Betreuenden bietet sich im Rahmen der gemeinsamen Erarbeitung der
Modelle die M\"oglichkeit und auch Herausforderung, speziell auf die individuelle
St\"arke unserer Gruppe einzugehen. So machten wir die Erfahrung, dass der
Motivationspegel und das Interesse besser gehalten werden kann, wenn wir auf
Formalismen wie etwa Summenzeichen -- das vielen Sch\"uler(inne)n noch nicht
vertraut ist -- zun\"achst verzichten. Zeigt sich, dass die Gruppe sehr stark
ist und das Optimierungsproblem vergleichsweise schnell gemeinsam formulieren
konnte, kann man sie anschlie\ss end zu einer kompakteren, anspruchsvolleren
Modellierung f\"uhren.
Unterschiede in der mathematischen St\"arke der Teilnehmenden, auch abh\"angig von
der Rahmenveranstaltung, sind nicht von der Hand zu weisen.  So wird mit dem
Girls'~Day oder auch der MINT\nobreakdash-Herbst\-uni\-ver\-si\-t\"at -- beide
Veranstaltungen dienen der Berufsorientierung -- ein allgemeineres Publikum
angesprochen, w\"ahrend der Tag der Mathematik sich speziell an
mathematikbegeisterte Sch\"uler(innen) wendet, die im gesamten Rahmen auch an
Mathematikwettbewerben teilnehmen.  Bez\"uglich der mathematischen F\"ahigkeiten
sehen wir keine geschlechtsspezifischen Unterschiede bei den Teilnehmenden.
Eine besondere Erfahrung f\"ur uns war das Optimierungslabor im Rahmen des
Hochschultages f\"ur im Fach Mathematik besonders begabte bayerische
Mittelstufensch\"uler(innen). Selbst die kompakte Formulierung der Modelle stellte f\"ur
die Teilnehmenden keine gr\"o\ss ere Schwierigkeit dar. Die Eigeninitiative war
au\ss ergew\"ohnlich hoch und es bestand sp\"urbares Interesse an weiterf\"uhrenden und
mathematisch sehr anspruchsvollen Aspekten. So fragte diese Gruppe immer wieder
danach, wie eine L\"osungsmethode f\"ur ILPs denn nun funktioniert (ein guter Grund,
Mathematik zu studieren!). Es wurde sogar (sinngem\"a\ss ) ein Strukturresultat
v\"ollig ungefragt gefunden \emph{und} dem Sinn nach korrekt nachgewiesen: wenn
alle St\"adte \emph{au\ss en} liegen (konvexe Lage), dann muss eine k\"urzeste Tour
au\ss en herum f\"uhren, da man sonst durch Austausch von sich kreuzenden Diagonalen
mit nicht kreuzenden Verbindungen die Tour verk\"urzen kann. Da waren wir {\glqq}platt{\grqq}!

\section{Abschlussbemerkungen}
\label{sec:Fazit}

Die Mathematik einmal ausschlie\ss lich als Sprache an der Schnittstelle
zwischen Mensch und Maschine zu betrachten, kann sich lohnen.
Existierende Mathematik-Software erm\"oglicht es, den Prozess der
Modellierung eines Problems vom Prozess der L\"osung zu trennen, um auch
ohne Kenntnisse der mathematischen L\"osungsverfahren zu einem
Endergebnis zu kommen.  Die Bedienung des Computers und der Software
kann dabei getrost dem Betreuungsteam \"uberlassen werden.

Unsere Erfahrungen zeigen n\"amlich, dass die eigenst\"andige,
handwerkliche Benutzung von Computern f\"ur eine Lerngruppe in einer
Veranstaltung mit Zeitrestriktionen auch zum Hemmschuh werden kann.
Es handelt sich wohl um ein weitverbreitetes Missverst\"andnis, dass die
F\"ahigkeit zur \emph{Benutzung} eines Computers \emph{die}
High-Tech-Kompetenz schlechthin darstellt.  Nach unserer Auffassung
stellt auch das Auffinden eines sinnvollen, konsistenten Wegs vom
Alltagsproblem zum Computer und von der Computerl\"osung zur\"uck zum
Ausgangsproblem eine kreative intellektuelle Herausforderung dar,
deren Bew\"altigung besonders nachhaltige Aha-Effekte hervorrufen
kann.

Diese Herausforderung l\"asst sich mit begreifbaren Materialien in
unserem Arbeitsgebiet lebendiger, kommunikativer und befriedigender
angehen als vor dem Bildschirm. Interessanterweise sind die dadurch
gef\"orderten und daf\"ur ben\"otigten F\"ahigkeiten gar nicht spezifisch f\"ur
Informationstechnologie. Logisches Denken, Initiative, Kritik- und
Urteilsf\"ahigkeit, Lese- und Zuh\"orverst\"andnis: das sind die
eigentlichen Erfolgsfaktoren f\"ur die Gruppen in unserem
Optimierungslabor.

Wir sind weder Lehrer(innen) noch Didaktiker(innen); unseren Konzepten und
Beobachtungen fehlt sicher eine einschl\"agige wissenschaftliche und
schulpraktische Fundierung. Deshalb wollen wir auch hier nicht besprechen,
inwieweit unser Projekt im Schulunterricht selbst anwendbar ist. Zu \"ahnlichen
Ideen gibt es Untersuchungen, z.\,B.{ }in \cite{Lutz-Westphal:Diss:2006}. Unser
Optimierungslabor ist letztlich ein zum Schulunterricht komplement\"ares
Schnupperangebot, das seinen universit\"aren und forschungsnahen Kontext auch
nicht verleugnen will. Trotzdem freuen wir uns \"uber Tipps und Anregungen der
Expert(inn)en!

Wir hoffen, dass unsere Erfahrungen zum Optimierungslabor als
Diskussionsbeitrag dienlich sind, und wir w\"urden uns freuen, die
Schulklassen von Leserinnen und Lesern vielleicht einmal bei uns
begr\"u\ss en zu d\"urfen.

\"Ubrigens: Die F\"unf-Damen-Konfiguration in Abbildung~\ref{fig:damen_index}
ist eine nicht-schlagende Damen-Konfiguration mit minimaler Anzahl von
Damen, so dass keine Dame mehr nicht-schlagend hinzugef\"ugt werden
kann.

\enlargethispage{2\baselineskip}
\footnotesize

\begin{thebibliography}{1}
  \providecommand{\bysame}{\leavevmode\hbox to3em{\hrulefill}\thinspace}
  \providecommand{\MR}{\relax\ifhmode\unskip\space\fi MR }
  \providecommand{\MRhref}[2]{%
    \href{http://www.ams.org/mathscinet-getitem?mr=#1}{#2}
  }
  \providecommand{\href}[2]{#2}
  \providebibliographyfont{name}{}%
  \providebibliographyfont{lastname}{}%
  \providebibliographyfont{title}{\emph}%
  \providebibliographyfont{jtitle}{\btxtitlefont}%
  \providebibliographyfont{etal}{}%
  \providebibliographyfont{journal}{}%
  \providebibliographyfont{volume}{\textbf}%
  \providebibliographyfont{ISBN}{\MakeUppercase}%
  \providebibliographyfont{ISSN}{\MakeUppercase}%
  \providebibliographyfont{url}{\url}%
  \providebibliographyfont{numeral}{}%
  \providecommand\btxprintamslanguage[1]{\ (#1)}
  \expandafter\btxselectlanguage\expandafter {\btxfallbacklanguage}

\expandafter\btxselectlanguage\expandafter {\btxfallbacklanguage}
\bibitem {bartlett+etal:sudoku}
\btxnamefont {Andrew\btxfnamespacelong C. \btxlastnamefont {Bartlett}},
  \btxnamefont {Timothy \btxlastnamefont {Chartier}}, \btxnamefont
  {Amy\btxfnamespacelong N. \btxlastnamefont {Langville}}\btxandcomma {}
  \btxandlong {} \btxnamefont {Timothy\btxfnamespacelong D. \btxlastnamefont
  {Rankin}}, \btxjtitlefont {\btxifchangecase {{An} {Integer} {Programming}
  {Model} for the {Sudoku} {Problem}}{{An} {Integer} {Programming} {Model} for
  the {Sudoku} {Problem}}}, \btxjournalfont {{Journal} of {Online Mathematics}
  and its {Applications}} \btxvolumefont {8} (2008).

\bibitem {0793.68113}
\btxnamefont {Axel \btxlastnamefont {Conrad}}, \btxnamefont {Tanja
  \btxlastnamefont {Hindrichs}}, \btxnamefont {Hussein \btxlastnamefont
  {Morsy}}\btxandcomma {} \btxandlong {} \btxnamefont {Ingo \btxlastnamefont
  {Wegener}}, \btxjtitlefont {\btxifchangecase {Solution of the knight's
  hamiltonian path problem on chessboards}{Solution of the knight's Hamiltonian
  path problem on chessboards}}, \btxjournalfont {Discrete Appl. Math.}
  \btxvolumefont {50} (1994), \btxnumbershort {.}~2, 125--134.

\bibitem {Kombinatorische_Optimierung_erleben}
\btxnamefont {Stephan\btxfnamespacelong (ed.) \btxlastnamefont {Hu{\ss}mann}}
  \btxandlong {} \btxnamefont {Brigitte\btxfnamespacelong (ed.)
  \btxlastnamefont {Lutz-Westphal}}, \btxtitlefont {\btxifchangecase
  {Kombinatorische {O}ptimierung erleben in {S}tudium und
  {U}nterricht}{Kombinatorische {O}ptimierung erleben in {S}tudium und
  {U}nterricht}}, \btxpublisherfont {Mathematik Erleben. Wiesbaden: Vieweg.
  xvi}, 2007.

\bibitem {koch:queens}
\btxnamefont {Thorsten \btxlastnamefont {Koch}}, \btxtitlefont
  {\btxifchangecase {Rapid mathematical programming}{Rapid Mathematical
  Programming}}, \btxphdthesis {}, Technische {Universit\"at} Berlin, 2004,
  {\latintext \btxurlfont{http://opus.kobv.de/zib/volltexte/2005/834/}},
  ZIB-Report 04-58.

\bibitem {koch:sudoku}
\bysame, \btxtitlefont {\btxifchangecase {Rapid mathematical programming or how
  to solve sudoku puzzles in a few seconds}{Rapid Mathematical Programming or
  How to Solve Sudoku Puzzles in a few Seconds}}, Operations Research
  Proceedings 2005\ (\btxnamefont {Hans\btxfnamespacelong Dietrich
  \btxlastnamefont {Haasis}}, \btxnamefont {Herbert \btxlastnamefont
  {Kopfer}}\btxandcomma {} \btxandlong {} \btxnamefont {J{\"o}rn
  \btxlastnamefont {Sch{\"o}nberger}}, \btxeditorsshort {.}), 2006, ZIB-Report
  05-51, \btxpagesshort {.}~21--26, \mbox{\btxISBN~\btxISBNfont
  {3-540-32537-9}}, {\latintext
  \btxurlfont{http://opus.kobv.de/zib/volltexte/2005/884/}}.

\bibitem {Lutz-Westphal:Diss:2006}
\btxnamefont {Brigitte \btxlastnamefont {Lutz-Westphal}}, \btxtitlefont
  {\btxifchangecase {{Kombinatorische} {Optimierung} -- {Inhalte} und
  {Methoden} f{\"u}r einen authentischen
  {Mathematikunterricht}}{{Kombinatorische} {Optimierung} -- {Inhalte} und
  {Methoden} f{\"u}r einen authentischen {Mathematikunterricht}}},
  \btxifchangecase {Dissertation}{Dissertation}, TU Berlin, 2006.

\bibitem {Rambau:Gelbe-Engel:2010}
\btxnamefont {J{\"o}rg \btxlastnamefont {Rambau}}, \btxtitlefont
  {\btxifchangecase {{Die} gelben {Engel} von {Noetham}}{{Die} gelben {Engel}
  von {Noetham}}}, {Besser} als {Mathe}: {Moderne} angewandte {Mathematik} aus
  dem {M\scshape atheon} zum {Mitmachen}\ (\btxnamefont {Katja \btxlastnamefont
  {Biermann}}, \btxnamefont {Martin \btxlastnamefont
  {Gr\"{o}tschel}}\btxandcomma {} \btxandlong {} \btxnamefont {Brigitte
  \btxlastnamefont {{Lutz-Westphal}}}, \btxeditorsshort {.}), \btxpublisherfont
  {Vieweg+Teubner}, 2010, \btxpagesshort {.}~59--74.

\bibitem {Rambau+Schwarz:ADAC-Aufzuege-Modelle:2008}
\btxnamefont {J{\"o}rg \btxlastnamefont {Rambau}} \btxandlong {} \btxnamefont
  {Cornelius \btxlastnamefont {Schwarz}}, \btxtitlefont {\btxifchangecase
  {Optimierte dynamische {Einsatzplanung} f{\"u}r {Gelbe} {Engel} und
  {Lastenaufz{\"u}ge}}{Optimierte dynamische {Einsatzplanung} f{\"u}r {Gelbe}
  {Engel} und {Lastenaufz{\"u}ge}}}, Die Kunst des Modellierens\ (\btxnamefont
  {Bernd \btxlastnamefont {Luderer}}, \btxeditorshort {.}), Teubner
  Studienb{\"u}cher Wirtschaftsmathematik, \btxpublisherfont {Vieweg+Teubner},
  2008, \btxpagesshort {.}~377--398.

\end{thebibliography}

\end{document}